\documentclass[12pt]{article}
\usepackage{latexsym}
\usepackage{amsmath}
\usepackage{amssymb}
\usepackage{amsfonts}
\usepackage{amsthm}
\usepackage[active]{srcltx}

\numberwithin{equation}{section}

\newtheorem{theorem}{Theorem}[section]
\newtheorem{lemma}[theorem]{Lemma}

\newtheorem{remark}[theorem]{Remark}
\newcommand{\eproof}{{\mbox{\ }~\hfill
\mbox{\large $\Box$} \par \vskip 10pt}}
\newcommand{\C}{{\mathbb C}}

\newcommand{\R}{{\mathbb R}}

\newcommand{\pf}{\noindent{\bf Proof}}

\title{Quantitative uniqueness estimates for the general second order elliptic equations}
\author{Ching-Lung Lin\thanks{Department of Mathematics, NCTS, National Cheng Kung University,
Tainan 701, Taiwan. Email:cllin2@mail.ncku.edu.tw}\qquad
Jenn-Nan Wang\thanks{Department of Mathematics, NCTS (Taipei), National Taiwan University,
Taipei 106, Taiwan. Email: jnwang@math.ntu.edu.tw}}

\date{}

\begin{document}
\maketitle

\begin{abstract}
In this paper we study quantitative uniqueness estimates of solutions to general second order elliptic equations with magnetic and electric potentials. We derive lower bounds of decay rate at infinity for any nontrivial solution under some general assumptions. The lower bounds depend on asymptotic behaviors of magnetic and electric potentials. The proof is carried out by the Carleman method and the bootstrapping arguments.
\end{abstract}

\section{Introduction}\label{sec1}
\setcounter{equation}{0}

In this paper we study the asymptotic behaviors of solutions to the general second order elliptic equation
\begin{equation}\label{1.1}
\begin{array}{l}
Pv+W(x)\cdot \nabla v+V(x)v+q(x)v=0\quad\text{in}\quad\Omega:=\R^n\setminus\bar{B},
\end{array}
\end{equation}
where $B$ is a bounded set in $\Omega$. Here $P(x,D)=\sum_{jk}a_{jk}(x)\partial_{j}\partial_{k}$ is uniformly
elliptic, i.e., for some $\lambda_0>0$
\begin{equation}\label{ellip}
\lambda_0|\xi|^2\le\sum_{jk}a_{jk}(x)\xi_j\xi_k\le\lambda_0^{-1}|\xi|^2\quad\forall\ x\in\Omega,\ \xi\in\R^n.
\end{equation}
and $a_{jk}(x)$ is Lipschitz continuous. We are interested in deriving lower bounds of the decay rate for any nontrivial solution $v$ to \eqref{1.1} under certain a priori assumptions. This kind of problem was originally posed by Landis in the 60's \cite{kl}. He conjectured that if $v$ is a bounded solution of
\begin{equation}\label{sch0}
\Delta v+q(x)v=0\quad\text{in}\quad \R^n
\end{equation}
with $\|q\|_{L^{\infty}}\le 1$ and $|v(x)|\le C\exp(-C|x|^{1+})$ for some constant $C$, then $v$ is identically zero. This conjecture was disproved by Meshkov \cite{me} who constructed a $q(x)$ and a nontrivial $v(x)$ with $|v(x)|\le C\exp(-C|x|^{4/3})$ satisfying \eqref{sch0}. He also proved that if $|v(x)|\le C_{k}\exp(-k|x|^{4/3})$ for all $k>0$ then $v\equiv 0$.  Note that $q(x)$ and $v(x)$ constructed by Meshkov are complex valued. In 2005, Bourgain and Kenig \cite{bou} derived a quantitative version of Meshkov's result in their resolution of Anderson localization for the Bernoulli model. Precisely, they showed that if $v$ is a bounded solution of $\Delta v+q(x)v=0$ in $\R^n$ satisfying $\|q\|_{L^{\infty}}\le 1$ and $v(0)=1$, then
\[
\inf_{|x_0|=R}\sup_{B(x_0,1)}|v(x)|\ge C\exp(-R^{4/3}\log R).
\]
In view of Meshkov's example, the exponent $4/3$ is optimal.

Recently, Davey \cite{de} derived similar quantitative asymptotic estimates for \eqref{1.1} with  $P=-\Delta$ and $q(x)=-E\in\C$, i.e.,
\begin{equation}\label{1.101}
-\Delta v+W(x)\cdot\nabla v+V(x)v=E v.
\end{equation}
To describe her result, we define
\begin{equation*}
I(x_0)=\int_{|y-x_0|<1}|v(y)|^2 dy
\end{equation*}
and
\begin{equation*}
M(t)=\inf_{|x_0|=t}I(x_0).
\end{equation*}
Assume that $|V(x)|\lesssim\langle x\rangle^{-N}$ and $|W(x)|\lesssim\langle x\rangle^{-P}$, where $\langle x\rangle=\sqrt{1+|x|^2}$. Then it was shown that for any nontrivial bounded solution $v$ of \eqref{1.101} with $v(0)=1$, we have
\begin{equation}\label{daest}
M(t)\gtrsim\exp(-Ct^{\beta_0}(\log t)^{b(t)}),
\end{equation}
where
\[
\beta_0=\max\{2-2P,\frac{4-2N}{3},1\}
\]
and $b(t)$ is either a constant $C$ or $C\log\log t$. Moreover, in \cite{de}, some Meshkov's type examples were constructed to ensure the optimality of \eqref{daest}. There are also some related qualitative results in \cite{cs}, \cite{efv}, \cite{ekpv11}, and \cite{fhho}. Especially, in \cite{cs} and \cite{fhho}, the authors studied the Schr\"odinger equation with potential $-\Delta v+V(x)v=E u$, where $|V(x)|\lesssim\langle x\rangle^{-N}$ with $0<N<1/2$ (in \cite{cs}) and $N\le 0$, $N>1/2$ (in \cite{fhho}). In addition to qualitative results, they also showed the optimality of $\beta_0$ (here $\beta_0=\max\{\frac{4-2N}{3},1\}$). For the case of $N=1/2$, the qualitative result was proved in \cite{ekpv11}.

In this work, we extend Davey's results to more general cases. Precisely, we consider the second order elliptic operator $P$ with more general assumptions on the asymptotic behaviors of $W$, $V$, and $q$. The main theorem is stated as follows.
\begin{theorem}\label{thm1.1}
Let $v\in H^1_{loc}(\Omega)$ be a nontrivial solution of
\eqref{1.1} satisfying
\begin{equation}\label{v}
|v(x)|\le \lambda|x|^{\alpha}\quad\text{with}\quad \alpha\ge 0
\end{equation}
for some $\lambda>0$. Assume that the ellipticity condition \eqref{ellip} holds and for $\epsilon>0$, $\kappa_1,\kappa_2,\kappa_3,\kappa_4\in\R$,
\begin{equation}\label{1.3}
\left\{
\begin{aligned}
&|W(x)|\leq \lambda |x|^{-\kappa_1},\,|V(x)|\leq \lambda |x|^{-\kappa_2},\\
&|\nabla a_{ij}(x)|\leq\lambda|x|^{-1-\epsilon},\\
&|q(x)|\le\lambda|x|^{\kappa_3},\, |\nabla q(x)|\leq \lambda|x|^{-\kappa_4}.
\end{aligned}\right.
\end{equation}
Denote
$\kappa_0=\max\{2-2\kappa_1,\frac{4-2\kappa_2}{3},\frac{2+\kappa_3}{2},\frac{3-\kappa_4}{2}\}$, $\kappa=\max\{\kappa_0,1\}$. Then we have that
\begin{itemize}
\item For $\kappa>1$ {\rm(}i.e., $\kappa_0>1${\rm)}, there exist $t_0$ depending on $\lambda_0$, $\lambda$,
$\epsilon$ and positive constants $C$, $C'$ such that
\begin{equation}\label{1.6}
M(t)\geq \exp\left(-Ct^{\kappa_0}(\log t)^{\gamma(t)}\right)\quad\text{for}\quad t\ge t_0
\end{equation}
with
\[
\gamma(t)=\frac{C'(\log t)(\log\log\log t)}{(\log\log t)^2},
\]
where $C=C(\lambda_0,\lambda,\epsilon,\kappa_1,\kappa_2,\kappa_3,\alpha)$ and $C'=C'(\lambda_0,\lambda,\epsilon,\kappa_1,\kappa_2,\kappa_3,\kappa_4,\alpha)$.

\item For $\kappa=1$ {\rm (}i.e., $\kappa_0\le 1${\rm)}, there exists a positive constant $C''$ such that
\begin{equation}\label{1.7}
M(t)\geq \exp\left(-C''t(\log t)^{\gamma(t)}\right)\quad\text{for}\quad t\ge  t_0,
\end{equation}
where $C''$ depends on $\alpha$, $\lambda_0$, $\lambda$, $\kappa_1, \kappa_2, \kappa_3$, $\epsilon$,
$$\left|\log\left(\min\{\inf_{t_0^{\frac{1}{1+\epsilon}}<|x|<t_0}\int_{|y-x|<1}|v(y)|^2dy,1\}\right)
\right|.$$ 
\end{itemize}
\end{theorem}

\begin{remark}
\begin{enumerate}

\item[{\rm 1.}] The condition on the decay rate of $\nabla a_{ij}$ was also used by T Nguyen {\rm\cite{ng}} in his proof of qualitative and quantitative Landis-Oleinik conjecture, the parabolic counterpart of Landis conjecture.

\item[{\rm 2.}] We have made general assumptions on the asymptotic behaviors of $W$, $V$, and $q$. They may grow in $|x|$. Our theorem provides quantitative uniqueness estimates for solutions of $-\Delta v+V(x) v=Ev$ with $|V|\lesssim |x|^m$, $m>0$.   Moreover, our method works for any $\kappa_0\in\R$ including the case $\kappa_0=1$ that is missing in {\rm\cite{de}}.

\end{enumerate}
\end{remark}

Similar to the arguments in \cite{bou} and \cite{de}, a key ingredient to prove Theorem~\ref{thm1.1} is the Carleman type estimate. Before applying the Carleman estimate, $v$ is shifted and rescaled appropriately. We will modify our ideas in \cite{lw}. Basically, we use the Carleman estimate for the shifted and rescaled solution on a ball of radius depending on $|x_0|$ (see the definition of $\Omega_t$ in Section~3). Note that in order to use the behaviors of coefficients in \eqref{1.3}, the radius of the ball is sufficiently small. In fact, after undoing sifting and rescaling, any point in $\Omega_t$ is at least $|x_0|^{1-\delta(x_0)}$ distance from the origin, where $\delta(x_0)=(\log\log |x_0|)^2/\log|x_0|$ (see Section~3). We then apply the Carleman estimate to derive three-ball inequalities in which we can estimate the $L^2$ bound of the solution in a unit ball centered at $x_0/|x_0|^{\delta(x_0)}$ by the $L^2$ bound of the solution in a unit ball centered at $x_0$ up to certain power (see \eqref{3.21}).  To obtain the desired estimates,  we want to apply bootstrapping arguments based on a chain of balls similar to what we did in \cite{lw}. An bootstrapping step was also used in \cite{de} to prove estimate \eqref{daest}. However, our method here is simpler than that of \cite{de}.

We now discuss the optimality of \eqref{1.6} and \eqref{1.7}, at least, for some simple cases.  It is readily seen that if $v(x)=\exp(-|x|^{1+\varepsilon})$ with $0<\varepsilon\le 1/2$, then $\Delta v+q(x)v=0$ in $\{|x|<1\}^c$ with $|q(x)|\sim |x|^{2\varepsilon}$ and $|\nabla q(x)|\sim|x|^{2\varepsilon-1}$. In this case, we can see that $\kappa=1+\varepsilon$. So the exponent $\kappa$ of \eqref{1.6} is optimal. This example also shows that if the first derivatives of the potential possess certain decaying property, we can break the $4/3$ barrier in the case of bounded potentials.  On the other hand, for $\varepsilon=0$, we obtain that $v=\exp(-|x|)$ satisfies $\Delta v+q(x)v=0$ in $\{|x|<1\}^c$ with $q(x)=-1+(n-1)|x|^{-1}$. Since we can write
\[
(\log t)^{\gamma(t)}=t^{\frac{C\log\log\log t}{\log\log t}},
\]
\eqref{1.7} is equivalent to
\[
M(t)\geq \exp\left(-Ct^{1+o(1))}\right).
\]
Thus, \eqref{1.7} is almost optimal.

This paper is organized as follows. In Section~2, we prove a Carleman estimate for the operator $P+q$, which plays an essential role in our proof. In Section~3, we begin to prove the main theorem, Theorem~\ref{thm1.1}, by deriving three-ball inequalities for solutions of \eqref{1.1}. In Section~4, we give detailed arguments of bootstrapping and complete the proof of Theorem~\ref{thm1.1}.

\section{Carleman estimate}\label{sec2}
\setcounter{equation}{0}

In this section, we would like to derive a Carleman estimate for $P+q$ with $q$ is $C^1$. Similar Carleman estimate for such operator with $P$ being the Laplace-Beltrami operator was also derived in \cite{ba} using Donnelly and Fefferman's approach \cite{dofe}.  Since we are working in the Euclidean space, we give a more elementary proof motived by the ideas in \cite{Reg1}. To begin, we introduce polar coordinates in ${\mathbb R}^n \backslash {\{0\}}$
by setting $x=r \omega$, with $r=|x|$,
$\omega=(\omega_1,\cdots,\omega_n)\in S^{n-1}$. Furthermore, using
new coordinate $t=\log r$, we can see that
$$\frac{\partial}{\partial x_j}=e^{-t}(\omega_j \partial_t +\Omega_j),\quad 1\le j\le n,$$
where $\Omega_j$ is a vector field in $S^{n-1}$. We could check that
the vector fields $\Omega_j$ satisfy
$$\sum_{j}\omega_j\Omega_j=0\quad\text{and}\quad\sum_{j}\Omega_j\omega_j=n-1.$$
Since $r\rightarrow 0$ iff $t\rightarrow {-\infty}$, we are mainly
interested in values of $t$ near $-\infty$.

It is easy to see that
\begin{equation*}
\frac{\partial ^2}{\partial x_j \partial x_{\ell}}=e^{-2t}(\omega_j
\partial_t -\omega_j +\Omega_j)(\omega_{\ell} \partial_t +\Omega_{\ell}),\quad 1\le j,\ell\le n.
\end{equation*}
and, therefore, the Laplacian becomes
\begin{equation}\label{2.1}
e^{2t}\Delta =\partial^2_t +(n-2)\partial_t +\Delta_\omega,
\end{equation}
where $\Delta_\omega=\Sigma_{j}\Omega^2_j$ denotes the
Laplace-Beltrami operator on $S^{n-1}$. We recall that the
eigenvalues of $-\Delta_\omega$ are $k(k+n-2), k\in \mathbb{N}$, and
the corresponding eigenspaces are $E_k$, where $E_k$ is the space of
spherical harmonics of degree $k$. It follows that
\begin{equation}\label{2.2}
\iint |\Delta_\omega v|^2 dt d\omega=\sum_{k\geq 0}k^2(k+n-2)^2
\iint | v_k |^2 dt d\omega
\end{equation}
and
\begin{equation}\label{2.3}
\sum_{j}\iint |\Omega_j v|^2 dt d\omega=\sum_{k\geq 0}k(k+n-2)\iint|
v_k |^2 dt d\omega,
\end{equation}
where $v_k$ is the projection of $v$ onto $E_k$.

Our aim is to derive Carleman-type
estimates with weights $\varphi_{\beta}=\varphi_{\beta}(x) =\exp
(-\beta\tilde{\psi}(x))$, where $\beta>0$ and $\tilde{\psi}(x)=\log
|x|+\log((\log |x|)^2)$.  For simplicity, we denote
$\psi(t)=t+\log t^2$, i.e., $\tilde{\psi}(x)=\psi(\log|x|)$.
\begin{lemma}\label{lem2.1}
Assume that $P=\sum_{jk}a_{jk}(x)\partial_j\partial_k$ is a second order elliptic operator satisfying
\[
\lambda_0|\xi|^2\le\sum_{jk}a_{jk}(x)\xi_j\xi_k\le\lambda_0^{-1}|\xi|^2\quad\forall\ x\in\R^n,\ \xi\in\R^n,
\]
and
\[
\|\nabla a_{jk}\|_{L^{\infty}(\R^n)}\le L.
\]
Then there exist a sufficiently small $r_1=r_1(\lambda_0,L)>0$ such that for all
$u\in U_{r_1}$ and $$\beta\geq (||q||_{L^\infty(\R^n)}+||\nabla
q||_{L^\infty(\R^n)})^{1/2},$$ we have
\begin{equation}\label{2.4}
\begin{aligned}
&\beta  \int \varphi^2_\beta (\log|x|)^{-2}|x|^{-n}(|x|^{2}|\nabla
u|^2+|u|^2)dx\\
&\le C_0\int \varphi^2_\beta|x|^{-n}|x|^{4}|Pu+q(x)u|^2dx,
\end{aligned}
\end{equation}
where $U_{r_1}=\{u\in C_0^{\infty}(\R^n\setminus\{0\}): \mbox{\rm
supp}(u)\subset B_{r_1}\}$ and $C_0=C_0(\lambda_0,L)>0$.
\end{lemma}

\pf. If we set $u=e^{\beta\psi(t)}v$ and $P_\beta
v=e^{-\beta\psi(t)} P (e^{\beta\psi(t)}v)$, then
\begin{eqnarray}\label{2.5}
&&e^{2t}P_\beta v+e^{2t}qv \notag\\
&=&\partial_t^2v+b\partial_tv+av+\Delta_\omega v+e^{2t}qv+\sum_{j+|\alpha|\le 2}C_{j\alpha}(t,\omega)\partial_t^j\Omega^\alpha v\notag \\
&&+\sum_{j+|\alpha|\le 1}C_{j\alpha}(t,\omega)\beta\psi'\partial_t^j\Omega^\alpha v+C_{20}(t,\omega)(\beta^2\psi'^2+\beta\psi'')v\notag\\
&=&\partial_t^2v+b\partial_tv+av+\Delta_\omega v+e^{2t}qv+Sv+hv,
\end{eqnarray}
where
\begin{equation*}
\left\{
\begin{aligned}
a&=\beta\psi''+\beta^2(\psi')^2+(n-2)\beta\psi'\\
&=(1+2t^{-1})^2\beta^2+(n-2)\beta+2(n-2)t^{-1}\beta-2t^{-2}\beta,\\
b&=2\beta\psi'+n-2\\
&=2\beta+4\beta t^{-1}+n-2,\\
S(v)&=\sum_{j+|\alpha|\le  2}C_{j\alpha}(t,\omega)\partial_t^j\Omega^\alpha v+C_{20}(t,\omega)(\beta^2\psi'^2+\beta\psi'')v,\\
h(v)&=\sum_{j+|\alpha|\le 1}C_{j\alpha}(t,\omega)\beta\psi'\partial_t^j\Omega^\alpha v,\\
C_{j\alpha}&=O(e^t), \partial_tC_{j\alpha}=O(e^t), \Omega^\alpha C_{j\alpha}=O(e^t).
\end{aligned}\right.
\end{equation*}
It is clear that \eqref{2.4} holds if for $t$ near $-\infty$ we have
\begin{align}\label{2.6}
&\iint |e^{2t}P_\beta v+e^{2t}qv|^2 dt d\omega\notag\\
&\geq C\left\{\iint \beta t^{-2}|\partial_tv|^2 dt d\omega+\beta\sum_{j}\iint t^{-2}|\Omega_jv|^2dtd\omega+\beta^3\iint t^{-2}|v|^2dtd\omega\right\}.
\end{align}

We obtain from \eqref{2.5} that
\begin{equation}\label{2.7}
|e^{2t}P_\beta v+e^{2t}qv|^2=|L(v)|^2+2b\partial_tvL(v)+2h(v)L(v)+|b\partial_tv+h(v)|^2,
\end{equation}
where $L(v):=\partial_t^2v+av+\Delta_\omega v+e^{2t}qv+S(v)$. Now we write
\begin{eqnarray}\label{2.8}
&& \iint 2b\partial_tvL(v) dt d\omega\notag\\
&=&\iint2b\partial_tv\partial_t^2vdtd\omega+\iint2abv\partial_tvdtd\omega+\iint2be^{2t}q\partial_tvvdtd\omega\notag\\
&&+\iint 2b\partial_tv\Delta_\omega vdtd\omega+\iint2b\partial_tvS(v)dtd\omega\notag.\\
\end{eqnarray}
By straightforward computations, we can get that
\begin{equation}\label{2.9}
\iint2b\partial_tv\partial_t^2vdtd\omega=-\iint(\partial_tb)|\partial_tv|^2dtd\omega,
\end{equation}
\begin{equation}\label{2.10}
\iint2abv\partial_tvdtd\omega=-\iint\partial_t(ab)|v|^2dtd\omega,
\end{equation}
\begin{equation}\label{2.11}
\iint2b\partial_tv\Delta_\omega vdtd\omega=\sum_{j}\iint(\partial_tb)|\Omega_j v|^2dtd\omega,
\end{equation}
and
\begin{equation}\label{2.12}
\iint2be^{2t}qv\partial_tvdtd\omega=-\iint(\partial_tb)e^{2t}q|v|^2dtd\omega-\iint\partial_t(e^{2t}q)b|v|^2dtd\omega.
\end{equation}
Combining \eqref{2.7} to \eqref{2.12} yields
\begin{eqnarray}\label{2.13}
&& \iint |e^{2t}P_\beta v+e^{2t}qv|^2 dt d\omega\notag\\
&=&\iint |L(v)|^2dt d\omega+4\beta\iint t^{-2}|\partial_tv|^2 dt d\omega\notag\\
&&+12\beta^3\iint t^{-2}(1+O(t^{-1}))|v|^2dtd\omega-4\beta\sum_{j}\iint t^{-2}|\Omega_jv|^2dtd\omega\notag\\
&&+4\beta\iint t^{-2}e^{2t}q|v|^2dtd\omega-\iint\partial_t(e^{2t}q)b|v|^2dtd\omega+\iint2b\partial_tvS(v)dtd\omega\notag\\
&&+\iint2h(v)L(v)dtd\omega+\iint|b\partial_tv+h(v)|^2dtd\omega.
\end{eqnarray}

Likewise, we write
\begin{equation}\label{2.14}
\begin{array}{l}
|L(v)|^2=|L(v)+3\beta t^{-2}v|^2-6\beta t^{-2}vL(v)-(3\beta t^{-2})^2|v|^2.
\end{array}
\end{equation}
It is easy to check that
\begin{align}\label{2.15}
&-6\beta\iint t^{-2}vL(v)dtd\omega\notag\\
=&-6\beta\iint t^{-2}v(\partial_t^2v+av+\Delta_\omega v+e^{2t}qv+S(v))dtd\omega\notag\\
=&-6\beta^3\iint t^{-2}(1+O(t^{-1}))|v|^2dt d\omega-6\beta\iint e^{2t}qt^{-2}|v|^2dt d\omega\notag\\
&+6\beta\iint t^{-2}|\partial_tv|^2 dt d\omega-12\beta\iint t^{-3}v\partial_tvdtd\omega+6\beta\sum_{j}\iint t^{-2}|\Omega_jv|^2dtd\omega\notag\\
&-6\beta\iint t^{-2}vS(v)dtd\omega.
\end{align}
From \eqref{2.13}-\eqref{2.15}, we have that for $t\le\tau$ ($\tau$ depends on $\lambda_0$, $L$)
\begin{align}\label{2.16}
& \iint |e^{2t}P_\beta v+e^{2t}qv|^2 dt d\omega\notag\\
\ge&\,\beta\iint t^{-2}|\partial_tv|^2 dt d\omega+2\beta\sum_{j}\iint t^{-2}|\Omega_jv|^2dtd\omega+3\beta^3\iint t^{-2}|v|^2dtd\omega\notag\\
&-2\beta\iint(q+\partial_tq)t^{-2}|v|^2dtd\omega+\iint2b\partial_tvS(v)dtd\omega\notag\\
&+\iint2h(v)L(v)dtd\omega-6\beta\iint t^{-2}vS(v) dtd\omega.
\end{align}
Using integration by parts and choosing an even smaller $\tau$, if necessary,  we can see that $\iint2b\partial_tvS(v)dtd\omega$, $\iint2h(v)L(v)dtd\omega$, and $6\beta\iint t^{-2}vS(v) dtd\omega$ are bounded by the first three terms on the right side of \eqref{2.16}. Therefore, taking
\[
\beta\geq\sqrt{\|q\|_{L^{\infty}}+\|\nabla q\|_{L^{\infty}}},
\]
\eqref{2.6} follows from \eqref{2.16}.
\eproof

\section{Proof of Theorem \ref{thm1.1} -- Part I: three-ball inequalities}\label{sec3}
\setcounter{equation}{0}

We begin to prove Theorem~\ref{thm1.1} in this section. As in \cite{bou}, \cite{de}, and \cite{lw}, the solution of \eqref{1.1} is shifted and rescaled properly. Fixing $x_0$ with  $|x_0|=t>>1$, we define
\begin{equation*}
\left\{
\begin{aligned}
&w(x)=u(atx+x_0),\,\tilde a_{jk}(x)=a_{jk}(atx+x_0),\,\tilde P(x,D)=\sum_{jk}\tilde a_{jk}(x)\partial_j\partial_k,\\
&\tilde W(x)=(at)W(atx+x_0),\,\tilde{V}(x)=(at)^2V(atx+x_0),\\
&\tilde{q}(x)=(at)^2q(atx+x_0),
\end{aligned}
\right.
\end{equation*}
where $a\geq
1/r_1$ will be determined later in the proof. Here $r_1$ is the constant appeared in Lemma~\ref{lem2.1}. We now denote
$$\Omega_{t}:=B_{\frac{1}{a}-\frac{1}{20at^{\delta}}}(0)=\{x: |x|<\frac{1}{a}-\frac{1}{20at^\delta}\}.$$  It follows from \eqref{1.1} that
\begin{equation}\label{3.1}
\tilde{P}w+\tilde{W}(x)\cdot \nabla w+\tilde{V}(x)w+\tilde{q}(x)w=0\quad\text{in}\quad\Omega_t.
\end{equation}
It is clear that $\tilde a_{jk}(x)$ satisfies \eqref{ellip} in $\Omega_t$ with same constant $\lambda_0$. Furthermore, in view of \eqref{1.3}, we have that
\begin{equation}\label{3.2}
\left\{
\begin{aligned}
&|\tilde{W}(x)|\leq 20\lambda at^{1-\kappa_1+\kappa_1\delta},\,|\tilde{V}(x)|\leq 20\lambda a^2t^{2-\kappa_2+\kappa_2\delta},\\
&|\nabla\tilde a_{ij}(x)|\leq 40\lambda at^{(1+\epsilon)\delta-\epsilon},\\
&|\tilde{q}(x)|\leq 2\lambda a^2t^{2+\kappa_3},\,|\nabla\tilde{q}(x)|\leq 20\lambda a^3t^{3-\kappa_4+\kappa_4\delta}.
\end{aligned}
\right.
\end{equation}
Unlike in \cite{lw}, where $\delta$ is a fixed constant, here we take $\delta=\delta(t)=\frac{(\log\log t)^2}{\log t}$. We now choose an $t_0$ such that $\log t_0\ge 1/r_1$ and
\begin{equation}\label{1.001}
-\epsilon_0:=(1+\epsilon)\delta(t_0)-\epsilon<0.
\end{equation}
By setting $a=t_0^{\epsilon_0}$, one can see that $at^{(1+\epsilon)\delta(t)-\epsilon}\le 1$ for all $t\ge t_0$. Let $r_1$ and $C_0$ be constants in Lemma~\ref{lem2.1} determined by $\lambda_0$ and $L=40\lambda$. Then the Carleman estimate \eqref{2.4} can be applied to $w$ in $\Omega_t$ for all $t\ge t_0$ with same $r_1$ and $C_0$. For simplicity, in this section, $C$ denotes a general constant whose value may vary from line to line. Furthermore, it depends on $\lambda_0$, $\lambda$, and $\epsilon$ unless indicated otherwise.

Besides of  the Carleman estimate \eqref{2.4}, we also need the following interior estimate for solutions of \eqref{3.1} in our proof.
\begin{lemma}\label{lem3.1}
For any $0<a_1<a_2$ such that $B_{a_2}\subset\Omega_{t}$ for $t>1$
and $a$ large enough, let $X=B_{a_2}\backslash \bar{B}_{a_1}$ and
$d(x)$ be the distant from $x\in X$ to $\mathbb{R}^n\backslash X$.
Then we have
\begin{eqnarray}\label{3.4}
&&(1+||\tilde{W}||^2_{L^\infty({X})})\int_{X}d(x)^{2}|\nabla w|^2dx\notag\\
&\le&
C\left(||\tilde{W}||^4_{L^\infty({X})}+||\tilde{V}||^2_{L^\infty({X})}+||\tilde{q}||^2_{L^\infty({X})}\right)\int_{X}|w|^2 dx
\end{eqnarray}
with $C=C(\lambda_0,L)$.
\end{lemma}
The lemma can be proved using similar arguments in \cite{luw2}. We omit the details here.

Now we are ready to apply \eqref{2.4} to $w$ solving \eqref{3.1}. Before doing so, we need to introduce a suitable cut-off function. Let $\chi(x)\in C^{\infty}_0
({\mathbb R}^n)$ satisfy $0\le\chi(x)\leq 1$ and
\begin{equation*}
\chi (x)=
\begin{cases}
\begin{array}{l}
0,\quad |x|\leq \frac{1}{8at},\\
1,\quad \frac{1}{4at}<|x|<\frac{1}{a}-\frac{3}{20at^\delta},\\
0,\quad |x|\geq \frac{1}{a}-\frac{2}{20at^\delta}.
\end{array}
\end{cases}
\end{equation*}
It is easy to see that for any multiindex $\alpha$
\begin{equation}\label{3.5}
\begin{cases}
|D^{\alpha}\chi|=O((at)^{|\alpha|})\quad \text{if}\quad \frac{1}{8at}\le |x|\le \frac{1}{4at},\\
|D^{\alpha}\chi|=O((at^\delta)^{|\alpha|})\quad \text{if}\quad \frac{1}{a}-\frac{3}{20at^\delta}\le |x|\le \frac{1}{a}-\frac{2}{20at^\delta}.
\end{cases}
\end{equation}
To use \eqref{2.4}, it suffices to take $\beta\geq\beta_1=\sqrt{20\lambda a^3}t^{\kappa_0}t^{|\kappa_4|\delta/2}$. Thus, we have
\begin{eqnarray}\label{3.6}
&&\int (\log|x|)^{-2}\varphi^2_\beta|x|^{-n}(\beta|x|^{2}|\nabla (\chi w)|^2+\beta^3|\chi w|^2)dx\notag\\
&\leq&C_0\int\varphi^2_\beta |x|^{-n}|x|^{4}|\tilde{P}(x) (\chi w)+\tilde{q}(\chi w)|^2 dx.
\end{eqnarray}
Using equation \eqref{3.1}, we obtain that
\begin{eqnarray}\label{3.7}
&&\int_{T}(\log|x|)^{-2}\varphi^2_\beta |x|^{-n}(\beta|x|^{2}|\nabla w|^2+\beta^3|w|^2)dx\notag\\
&\leq& \int\varphi^2_\beta (\log|x|)^{-2}|x|^{-n}(\beta|x|^{2}\nabla (\chi w)|^2+\beta^3|\chi w|^2)dx\notag\\
&\leq &C_0\int\varphi^2_\beta |x|^{-n}|x|^{4}|\chi (\tilde{W}(x)\cdot \nabla w+\tilde{V}(x)w)|^2 dx\notag\\
&&+C_0\int \varphi^2_\beta|x|^{4-n}|[\tilde{P},\chi ]w|^2,\notag\\
\end{eqnarray}
where $T$ denotes the domain
$\{x:\frac{1}{4at}<|x|<\frac{1}{a}-\frac{3}{20at^\delta}\}$. To
simplify the notations, we denote $Y=\{x:\frac{1}{8at}\le |x|\le
\frac{1}{4at}\}$ and $Z=\{x:\frac{1}{a}-\frac{3}{20at^\delta}\le
|x|\le \frac{1}{a}-\frac{2}{20at^\delta}\}$. By \eqref{3.2} and
estimates \eqref{3.5}, we deduce from \eqref{3.7} that
\begin{eqnarray}\label{3.8}
&& \int_{T}(\log|x|)^{-2}\varphi^2_\beta |x|^{-n}(\beta|x|^{2}|\nabla w|^2+\beta^3|w|^2)dx\notag\\
&\leq& C'\int_{T}\varphi^2_\beta |x|^{-n}|x|^{4}(a^2t^{2-2\kappa_1+2\kappa_1\delta}|\nabla w|^2+a^4t^{4-2\kappa_2+2\kappa_2\delta}|w|^2) dx\notag\\
&&+C\int_{Y}\varphi^2_\beta |x|^{-n}|\tilde{U}|^2 dx+C\|\tilde W\|_{L^{\infty}}^2\int_{Z}\varphi^2_\beta |x|^{-n}|x|^{2}|\nabla w|^2 dx\notag\\
&&+C\|\tilde V\|_{L^{\infty}}^2\int_{Z}\varphi^2_\beta |x|^{-n}|w|^2 dx,
\end{eqnarray}
where $|\tilde{U}(x)|^2=|x|^{2}|\nabla w|^2+|w|^2$, $C'=C'(\lambda_0,\lambda)$, $C=C(\lambda_0,\lambda,a)$. From now on, $\|\cdot\|_{L^{\infty}}$ is taken over $\Omega_t$.

Taking a larger $t_0$ (recall $a=t_0^{\epsilon_0}$), if necessary, we can obtain that $|x|^{2}(\log|x|)^{2}C'\leq \frac 12$ for all $x\in T$. Additionally, we choose
$\beta\geq\beta_2:= a^2t^{\kappa_0+\kappa_s\delta}$, where $\kappa_s=\max\{2|\kappa_1|,2|\kappa_2|/3,|\kappa_4|/2\}$. then the first term on the right hand side
of \eqref{3.8} can be absorbed by the left hand side of \eqref{3.8}.
With the choices described above, we obtain from \eqref{3.8} that
\begin{eqnarray}\label{3.9}
&&\beta^3(b_1)^{-n}(\log b_1)^{-2}\varphi^2_\beta(b_1)\int_{\frac{1}{2at}<|x|<b_1}|w|^2dx\notag\\
&\leq &\beta^3\int_{T}(\log|x|)^{-2}\varphi^2_\beta |x|^{-n}|w|^2dx\notag\\
&\leq &Cb_2^{-n}\varphi^2_\beta(b_2)\int_{Y}|\tilde{U}|^2 dx+C\|\tilde W\|_{L^{\infty}}^2b_3^{-n}\varphi^2_\beta(b_3)\int_{Z}|x|^{2}|\nabla w|^2 dx\notag\\
&&+C\|\tilde V\|_{L^{\infty}}^2b_3^{-n}\varphi^2_\beta(b_3)\int_{Z}|w|^2 dx,
\end{eqnarray}
where $b_1=\frac{1}{a}-\frac{8}{20at^\delta}$, $b_2=\frac{1}{8at}$
and $b_3=\frac{1}{a}-\frac{3}{20at^\delta}$.

Using \eqref{3.4}, we can control $|\tilde{U}|^2$ terms on the
right hand side of \eqref{3.5}. Indeed, let
$X=Y_1:=\{x:\frac{1}{16at}\le |x|\le \frac{1}{2at}\}$, then we can
see that
$$d(x)\ge C|x|\quad\text{for all}\quad x\in Y,$$ where $C$ an absolute constant. Therefore, \eqref{3.4} implies
\begin{eqnarray}\label{3.10}
&&\int_{Y}|x|^{2}|\nabla w|^2dx\notag\\
&\le&C\int_{Y_1}d(x)^{2}|\nabla w|^2dx\notag\\
&\le&C(\|\tilde{W}\|^4_{L^\infty}+\|\tilde{V}\|^2_{L^\infty}+\|\tilde{q}\|^2_{L^\infty})\int_{Y_1}|w|^2 dx,
\end{eqnarray}
where $C=C(\lambda_0,\lambda,a)$. On the other hand, let
$X=Z_1:=\{x:\frac{1}{2a}\le |x|\le
\frac{1}{a}-\frac{1}{20at^\delta}\}$, then
$$d(x)\ge Ct^{-\delta}|x|\quad\text{for all}\quad x\in Z,$$ where $C$ another absolute constant. Thus, it follows from \eqref{3.4} that
\begin{eqnarray}\label{3.11}
&&\int_{Z}|x|^{2}|\nabla w|^2dx\notag\\
&\le&Ct^{2\delta}\int_{Z_1}d(x)^{2}|w|^2dx\notag\\
&\le&Ct^{2\delta}(\|\tilde{W}\|^4_{L^\infty}+\|\tilde{V}\|^2_{L^\infty}+\|\tilde{q}\|^2_{L^\infty})\int_{Z_1}|w|^2dx.
\end{eqnarray}
Here, $C$ also depends on $\lambda_0$, $\lambda$. Combining \eqref{3.9}, \eqref{3.10}, and \eqref{3.11} leads to
\begin{eqnarray}\label{3.12}
&&b_1^{-2\beta-n}(\log b_1)^{-4\beta-2}\int_{\frac{1}{2at}<|x|<b_1}|w|^2dx\notag\\
&\leq &C''t^{p_1}b_2^{-n}\varphi^2_\beta(b_2)\int_{Y_1}|w|^2 dx+C''t^{p_2}b_3^{-n}\varphi^2_\beta(b_3)\int_{Z_1}|w|^2 dx,
\end{eqnarray}
where $C''=C''(\lambda_0,\lambda,a)$, $p_1=p_1(\kappa_1,\kappa_2,\kappa_3)$, $p_2=p_2(\kappa_1,\kappa_2,\kappa_3)$. Notice that \eqref{3.12} holds for all $\beta\ge\beta_2$.

Changing $2\beta+n$ to $\beta$, \eqref{3.12} becomes
\begin{eqnarray}\label{3.13}
&&b_1^{-\beta}(\log b_1)^{-2\beta+2n-2}\int_{\frac{1}{2at}<|x|<b_1}|w|^2dx\notag\\
&\leq & C''t^{p_1}b_2^{-\beta}(\log b_2)^{-2\beta+2n}\int_{Y_1}|w|^2 dx\notag\\
&&+C''t^{p_2}b_3^{-\beta}(\log b_3)^{-2\beta+2n}\int_{Z_1}|w|^2 dx.
\end{eqnarray}
Recall that $\delta=\delta(t)=\frac{(\log\log t)^2}{\log t}$. By taking $t_0$ sufficiently large, if necessary, we can see that for $t\ge t_0$
\begin{equation}\label{3.14}
\begin{cases}
\frac{1}{a}-\frac{1}{a{t}^\delta}+\frac{1}{a{t}}\leq \frac{1}{a}-\frac{8}{20a{t}^\delta},\\
\frac{1}{a}-\frac{1}{a{t}^\delta}-\frac{1}{a{t}}\geq \frac{1}{2a{t}}.
\end{cases}
\end{equation}
In view of \eqref{3.14}, dividing $b_1^{-\beta}(\log b_1)^{-2\beta+2n-2}$ on the both sides
of \eqref{3.13} and noting that we can let $\beta_2\ge n-1$, i.e.,
$2\beta-2n+2>0$ for all $\beta\ge\beta_2$, we obtain that
\begin{eqnarray}\label{3.15}
&&\int_{|x+\frac{b_4x_0}{t}|<\frac{1}{at}}|w(x)|^2dx\notag\\
&\leq &\int_{\frac{1}{2at}<|x|<b_1}|w(x)|^2dx\notag\\
&\leq & C''t^{p_1}(\log b_1)^{2}(b_1/b_2)^{\beta}\int_{Y_1}|w|^2dx\notag\\
&&+C''t^{p_2}(b_1/b_3)^\beta (\log b_3)^2[\log b_1/\log b_3]^{2\beta-2n+2}\int_{Z_1}|w|^2 dx\notag\\
&\leq & C''t^{p_1}(\log b_1)^{2}(8t)^{\beta}\int_{|x|<\frac{1}{at}}|w(x)|^2dx\notag\\
&&+C''t^{p_2}(\log b_3)^2(b_1/b_5)^{\beta}\int_{Z_1}|w(x)|^2 dx,
\end{eqnarray}
where $b_4=\frac{1}{a}-\frac{1}{at^\delta}$ and
$b_5=\frac{1}{a}-\frac{6}{20at^\delta}$. In deriving the third
inequality above, we use the fact that if $a(=t_0^{\epsilon_0})$ is sufficiently large, then
$$
0\le(\frac{b_5}{b_3})(\frac{\log b_1}{\log
b_3})^2\le 1
$$
for all $t\ge t_0$. From now on we fix $a$ and hence $t_0$. It is helpful to remind that $t_0$ depends on $\lambda_0$ and $\lambda$.  Having fixed constant $a$, $|\log
b_1|$ and $|\log b_3|$ can be bounded by a positive constant. Thus, \eqref{3.15} is
reduced to
\begin{eqnarray}\label{3.16}
\int_{|x+\frac{b_4x_0}{t}|<\frac{1}{at}}|w(x)|^2dx&\leq& Ct^{p_1}(8t)^{\beta}\int_{|x|<\frac{1}{at}}|w(x)|^2dx\notag\\
&&+Ct^{p_2}(b_1/b_5)^{\beta}\int_{Z_1}|w(x)|^2 dx,
\end{eqnarray}
where $C=C(\lambda_0,\lambda,\epsilon)$.

Using \eqref{3.16}, \eqref{v}, rescaling $w$ back to $v$,  and replacing $\beta$ by $\beta-p_1$ (that is, taking $\beta\ge\beta_2+p_1$), we have that
\begin{equation}\label{3.17}
I(t^{1-\delta}y_0)\leq C(8t)^{\beta}I(ty_0)+Ct^{p}\left(\frac{t^\delta}{t^\delta+0.1}\right)^{\beta},
\end{equation}
where $x_0=ty_0$, $C=C(\lambda_0, \lambda, \kappa_1, \kappa_2, \kappa_3, \epsilon)$, and $p=p(\kappa_1,\kappa_2,\kappa_3,\alpha)$.  For simplicity, by
denoting
\begin{equation*}
A(t)=\log 8t,\quad B(t)=\log \left(\frac{t^\delta+0.1}{t^\delta}\right),
\end{equation*}
\eqref{3.17} becomes
\begin{equation}\label{3.18}
I(t^{1-\delta}y_0)\leq C\Big{\{}\exp(\beta A(t))I(ty_0)+t^p\exp(-\beta B(t))\Big{\}}.
\end{equation}

Now, we consider two cases. If
$$ \exp({\beta}_2 A(t))I(ty_0)\geq t^p\exp(-{\beta}_2 B(t)),$$
then we have
\begin{equation*}
I(x_0)=I(ty_0)\ge t^p\exp(-\beta_2
(A(t)+B(t)))=t^p(8t)^{-\beta_2}\left(\frac{t^{\delta}+0.1}{t^{\delta}}\right)^{-\beta_2},
\end{equation*}
that is
\begin{equation}\label{3.19}
I(ty_0)\ge t^{-C\beta_2}=t^{-Ct^{\kappa_0+\kappa_s\delta}}\ge\exp(-Ct^{\kappa_0+\kappa_s\delta}\log t)
\end{equation}
for all $t\ge t_0$, where $C=C(\lambda_0,\lambda,\epsilon,\kappa_1,\kappa_2,\kappa_3,\alpha)$.  On the other hand, if
$$\exp(\beta_2 A(t))I(ty_0)< t^p\exp(-{\beta}_2B(t)),$$
then we can pick a $\tilde\beta>\beta_2$ such that
\begin{equation}\label{3.20}
\exp(\tilde\beta A(t))I(ty_0)=t^p\exp(-\tilde\beta B(t)).
\end{equation}
Solving $\tilde\beta$ from \eqref{3.20} and using \eqref{3.18}, we
have that
\begin{eqnarray}\label{3.21}
I(t^{1-\delta}y_0)&\leq& C\exp(\tilde\beta A(t))I(ty_0)\notag\\
&=& C\left(I(ty_0)\right)^{\tau}(t^{p})^{1-\tau}\notag\\
&\leq& Ct^{p}\left(I(ty_0)\right)^{\tau},
\end{eqnarray}
where $\tau=\frac{B(t)}{A(t)+B(t)}$. This estimate will serve as a building block in the bootstrapping step in the next section.

\section{Proof of Theorem~\ref{thm1.1} -- Part II: bootstrapping}\label{sec4}

In the previous section, we see that \eqref{3.19} gives us the desired estimate. However, we need to work harder to derive the wanted estimate from \eqref{3.21}. We first observe that for $t^{\frac{2}{(\log\log t)^2}}\le \hat{t}\le t$ we have.
\begin{equation*}
\begin{cases}
\frac{1}{a}-\frac{1}{a{\hat{t}}^\delta}+\frac{1}{a{\hat{t}}}\leq \frac{1}{a}-\frac{8}{20a{\hat{t}}^\delta},\\
\frac{1}{a}-\frac{1}{a{\hat{t}}^\delta}-\frac{1}{a{\hat{t}}}\geq \frac{1}{2a{\hat{t}}}.
\end{cases}
\end{equation*}
We now let $|x_0|=t$ with
\begin{equation}\label{4.001}
t_0\le t^{\frac{2}{(\log\log t)^2}},
\end{equation}
then we can write
\begin{equation}\label{3.22}
t=\mu^{\left((1-\delta)^{-s}\right)}
\end{equation}
for some positive integer $s$ and
\begin{equation}\label{4.002}
(t^{\frac{2}{(\log\log t)^2}})^{1-\delta}\leq\mu\leq
t^{\frac{2}{(\log\log t)^2}}.
\end{equation}
For simplicity, we define
$d_j=\mu^{\left((1-\delta)^{-j}\right)}$ and
$\tau_j=\frac{B(d_j)}{A(d_j)+B(d_j)}$ for $j=1,2\cdots s$. Define
$$\tilde{J}=\{1\leq j\leq s : \exp(a^2d_j^{\kappa_0+\kappa_s\delta} A(d_j))I(d_{j}y_0)\geq (d_j)^p\exp(-a^2d_j^{\kappa_0+\kappa_s\delta} B(d_j))\},$$ where $y_0=x_0/t$ as before. Note that $a$ is a fixed constant depending on $\lambda_0$, $\lambda$, $\epsilon$. Now, we divide it into two cases. If $\tilde{J}=\emptyset$, we only
need to consider \eqref{3.21}. Using \eqref{3.21} iteratively
starting from $t=d_1$, we obtain that
\begin{eqnarray}\label{3.23}
I(\mu y_0)&\leq& C(d_1^{p})\left(I(d_1y_0)\right)^{\tau_1}\notag\\
&\leq& C^s(d_1d_2\cdots d_s)^{p}\left(I(x_0)\right)^{\tau_1\tau_2\cdots\tau_s}.
\end{eqnarray}
It is easy to check that $s\leq C\log t(\log\log\log t/(\log\log t)^2)$ for some absolute constant $C$.
From \eqref{3.23} we have that
\begin{eqnarray}\label{3.24}
I(\mu y_0)&\leq& {C}^{s}t^{p/\delta}\left(I(x_0)\right)^{\tau_1\tau_2\cdots\tau_s}\notag\\
&\le&t^{C\log t/(\log\log t)^2}\left(I(x_0)\right)^{\tau_1\tau_2\cdots\tau_s}.
\end{eqnarray}
Hereafter, $C=C(\lambda_0, \lambda, \epsilon, \kappa_1,\kappa_2,\kappa_3, \alpha)$, unless indicated otherwise. We now estimate
\[
\frac{1}{\tau_j}=\frac{\log(8d_j)+\log(1+0.1d_j^{-\delta})}{\log(1+0.1d_j^{-\delta})}\le \frac{2\log(8d_j)}{\log(1+0.1d_j^{-\delta})}
\le 40d_j^{\delta}\log(d_j).
\]
and thus
\begin{eqnarray}\label{3.25}
\frac{1}{\tau_1\tau_2\cdots\tau_s}&\le& 40^s(\log t)^s(d_1\cdots d_s)^{\delta}\notag\\
{}&\le&t\omega(t),
\end{eqnarray}
where 
\begin{equation}\label{om}
\omega(t)=t^{\frac{C\log\log\log t}{\log\log t}}. 
\end{equation}
Raising both sides of
\eqref{3.24} to the power $\frac{1}{\tau_1\tau_2\cdots\tau_s}$ and
using \eqref{3.25}, we obtain that
\begin{equation*}
\left(\min\{I(\mu y_0),1\}\right)^{t\omega(t)}\leq t^{(C\log t/(\log\log t)^2)t\omega(t)}I(x_0),
\end{equation*}
i.e.,
\begin{equation}\label{3.26}
I(x_0)\ge \exp({-Ct\omega(t)})\left(\min\{I(\mu y_0),1\}\right)^{t\omega(t)}.
\end{equation}

Next, if $\tilde{J}\neq\emptyset$, let $l$ be the largest integer in
$J$. Estimate \eqref{3.19} implies
\begin{equation}\label{3.27}
I(d_ly_0)\ge\exp(-Cd_l^{\kappa_0+\kappa_s\delta}\log d_l)\ge\exp(-Ct^{\kappa_0+\kappa_s\delta}\log t).
\end{equation}
Note that here $\delta=\delta(t)$. As in \eqref{3.23}, iterating \eqref{3.21} starting from $t=d_{l+1}$ yields
\begin{eqnarray}\label{3.28}
I(d_ly_0)&\le& C^{s-l}(d_{l+1}\cdots
d_s)^{p}\left(I(x_0)\right)^{\tau_{l+1}\cdots\tau_s}\notag\\
&\le&{C}^{s}(t/d_l)^{p/\delta}\left(I(x_0)\right)^{\tau_{l+1}\cdots\tau_s}\notag\\
&\le&t^{C\log t/(\log\log t)^2}\left(I(x_0)\right)^{\tau_{l+1}\cdots\tau_s}.
\end{eqnarray}
It is enough to assume $I(d_ly_0)<1$. Repeating the computations in
\eqref{3.25}, we can see that
\begin{equation}\label{3.29}
\frac{1}{\tau_{l+1}\cdots\tau_s}\le(t/d_l)\omega(t).
\end{equation}
Hence, combining \eqref{3.27}, \eqref{3.28} and using \eqref{3.29},
we get that
\begin{eqnarray}\label{3.30}
\exp(-Ct^{\kappa_0}{\omega}(t))\leq \left(I(x_0)\right).
\end{eqnarray}
Here $\omega(t)$ is given as in \eqref{om}, but with $C=C(\lambda_0, \lambda, \epsilon, \kappa_1,\kappa_2,\kappa_3, \kappa_4, \alpha)$.

The last estimate \eqref{3.30} gives us the desired bound. We now focus on \eqref{3.26}.  In view of \eqref{4.002}, if $\mu$ satisfies
\[
(t^{\frac{2}{(\log\log t)^2}})^{1-\delta}\le\mu\le t_0,
\]
then we are done. Note that $t_0^{\frac{1}{1+\epsilon}}\le(t^{\frac{2}{(\log\log t)^2}})^{1-\delta}$ due to \eqref{1.001} and \eqref{4.001}. So we now consider $\mu>t_0$. We need another bootstrapping argument. Let us first rename $\tilde{t}_{0}=t$, $\delta_1=\delta(t)$, $s_1=s$ (as in \eqref{3.22}). We then denote $\tilde{t}_1=\mu$ and
$\delta_2=\delta_2(\tilde{t}_1)=\frac{(\log\log \tilde{t}_1)^2}{\log
\tilde{t}_1}$. As before, we write
\begin{equation*}
\tilde{t}_1=\tilde{t}_2^{\left((1-\delta_2)^{-s_2}\right)}
\end{equation*}
for some positive integer $s_2$ and $(\tilde{t}_1^{\frac{2}{(\log\log \tilde{t}_1)^2}})^{1-\delta_2}\leq\tilde{t}_2\leq \tilde{t}_1^{\frac{2}{(\log\log \tilde{t}_1)^2}} $.
Inductively, we denote
$\delta_k=\delta_k(\tilde{t}_{k-1})=\frac{(\log\log
(\tilde{t}_{k-1}))^2}{\log (\tilde{t}_{k-1})}$ and  write
\begin{equation*}
\tilde{t}_{k-1}=\tilde{t}_k^{\left((1-\delta_k)^{-s_k}\right)}
\end{equation*}
with positive constant $s_k$ and $(\tilde{t}_{k-1}^{\frac{2}{(\log\log \tilde{t}_{k-1})^2}})^{1-\delta_k}\leq\tilde{t}_k\leq \tilde{t}_{k-1}^{\frac{2}{(\log\log \tilde{t}_{k-1})^2}}$ for $k=1,2,\cdots$. It is easily seen that there exists an $m$ such that
\[
t_0^{\frac{1}{1+\epsilon}}\le \tilde t_m\le t_0.
\]
Indeed, we have $m\le\log t$.

Now we are ready to perform bootstrapping using either \eqref{3.26} or \eqref{3.30}. It is enough to treat the case where we have \eqref{3.26} all the way until $\tilde t_m$, namely,
\begin{equation}\label{4.010}
\begin{aligned}
&I(x_0)\\
&\ge e^{-Ct\omega(t)}e^{-C\tilde t_1\omega(\tilde t_1)t\omega(t)}e^{-C\tilde t_2\omega(t_2)\tilde t_1\omega(\tilde t_1)t\omega(t)}\cdots e^{-C\tilde t_{m-1}\omega(\tilde t_{m-1})\cdots\tilde t_2\omega(t_2)\tilde t_1\omega(\tilde t_1)t\omega(t)}\\
&\quad\times\left(\min\{I(\tilde t_m y_0),1\}\right)^{\tilde t_{m-1}\omega(\tilde t_{m-1})\cdots\tilde t_2\omega(t_2)\tilde t_1\omega(\tilde t_1)t\omega(t)}.
\end{aligned}
\end{equation}
In view of
\[
\omega(\tilde t_k)\le\tilde t_k,\quad\tilde t_{k+1}^2\le\tilde t_k\quad\text{for}\ 1\le k\le(m-1),
\]
we deduce that
\[
\tilde t_{m-1}\omega(\tilde t_{m-1})\cdots\tilde t_2\omega(t_2)\tilde t_1\omega(\tilde t_1)t\omega(t)\le\tilde t_1^4t\omega(t).
\]
Multiplying all terms in \eqref{4.010} implies
\begin{equation}\label{4.011}
I(x_0)\ge \exp(-Ct{\omega}(t)),
\end{equation}
where $C=C(\lambda_0, \lambda, \epsilon, \kappa_1,\kappa_2,\kappa_3, \alpha)$ and
$$\left|\log\left(\min\{\inf_{t_0^{\frac{1}{1+\epsilon}}<|x|<t_0}\int_{|y-x|<1}|v(y)|^2dy,1\}\right)
\right|$$
and $\omega(t)$ is defined as in \eqref{om} with $C$ depending on the same parameters. 

On the other hand, we stop the bootstrapping process whenever \eqref{3.30} is satisfied. Similar computations give the following bound
\begin{equation}\label{4.020}
I(x_0)\ge\exp(-C't{\omega}(t)),
\end{equation}
where $C'=C'(\lambda_0,\lambda,\epsilon,\kappa_1,\kappa_2,\kappa_3,\alpha)$ and the constant $C$ in $\omega(t)$ depends on $\lambda_0,\lambda,\epsilon,\kappa_1,\kappa_2,\kappa_3,\kappa_4,\alpha$. Notice that $$\omega(t)=t^{\frac{C\log\log\log t}{\log\log t}}=(\log t)^{\frac{C(\log t)(\log\log\log t)}{(\log\log t)^2}}.$$  Therefore, \eqref{3.30} gives the estimate for $\kappa>1$ and \eqref{3.30}, \eqref{4.011}, \eqref{4.020} lead to the estimate for $\kappa=1$.

\section*{Acknowledgements}
CL was supported in part by the National Science Council of Taiwan. JW was supported partially by the National Science Council of Taiwan grant 99-2115-M-002-006-MY3.

\end{document}